\documentclass[reqno,11pt]{amsproc}
\usepackage{amssymb}
\usepackage{fancybox}
\usepackage{array}
\usepackage{graphicx}
\usepackage{amscd}
\usepackage{amsmath}

\theoremstyle{plain}

\newtheorem{lemma}{Lemma}

\newtheorem{proposition}{Proposition}

\numberwithin{equation}{section}

\def\supp{\operatorname{supp}}

\def\OB{\operatorname{\O}}

\def\Nset{\mbox{I\kern-.21em N}}
\def\RE{{\mbox{\rm I\kern-.21em R}}}
\def\ZZ{{\mbox{\sf Z\kern-.45em Z}}}
\def\vv{\kern.344em{\rule[.18ex]{.075em}{1.32ex}}\kern-.344em}

 \def\O{\mathcal O}
  
\def\<{\langle} \def\>{\rangle}

\begin{document}

\title{One dimensional inverse problem in photoacoustic. Numerical testing.}
\author{D. Langemann}%
\address{Technische Universit{\"a}t Braunschweig, Inst. Computational
    Mathematics, AG PDE, Universit{\"a}tsplatz 2, 38106 Braunschweig, Germany}
\email{ d.langemann@tu-bs.de  }%
\author{A.S. Mikhaylov} 
\address{St. Petersburg   Department   of   V.A. Steklov    Institute   of   Mathematics
    of   the   Russian   Academy   of   Sciences, 7, Fontanka, 191023
    St. Petersburg, Russia and Saint Petersburg State University,
    St.Petersburg State University, 7/9 Universitetskaya nab., St.
    Petersburg, 199034 Russia.} \email{mikhaylov@pdmi.ras.ru}
\author{V.S. Mikhaylov} 
\address{St.Petersburg   Department   of   V.A.Steklov    Institute   of   Mathematics
    of   the   Russian   Academy   of   Sciences, 7, Fontanka, 191023
    St. Petersburg, Russia and Saint Petersburg State University,
    St.Petersburg State University, 7/9 Universitetskaya nab., St.
    Petersburg, 199034 Russia.} \email{ vsmikhaylov@pdmi.ras.ru }

\date{May, 2018}

\maketitle

\hfill {\bf Dedicated to the memory of A.P. Kachalov.}

\bigskip

{\bf Abstract.} We consider the problem of reconstruction of
Cauchy data for the wave equation in $\mathbb{R}^1$  by the
measurements of its solution on the boundary of the finite
interval. This is a one-dimensional model for the multidimensional
problem of photoacoustics, which was studied in \cite{BLMM}. We
adapt and simplify the method for one-dimensional situation and
provide the results on numerical testing to see the rate of
convergence and stability of the procedure. We also give some
hints on how the procedure of reconstruction can be simplified in
2d and 3d cases.

\section{Introduction}
Photoacoustic  is a method of visualization for obtaining optical
images in a scattering medium. Ultrashort pulses of laser
radiation causes thermoelastic stresses in the light absorption
region. Such thermal expansion causes the propagation of
ultrasonic waves in the medium, which can be registered and
measured. After that, using the measurements, mathematical
algorithms and computers one can produce an image, in this case
the technology called Photoacoustic tomography.

In what follows we assume that the visualized medium is homogeneous with respect
to ultrasound, put the speed of sound equal to 1 and write down the
standard wave equation for a pressure wave
\begin{equation}
    \label{wave_eq}
    \left\{\begin{array}l u_{tt}-\Delta u=0,\quad x\in
        \mathbb{R}^n,\,\,t>0,\\
        u\bigl|_{t=0}=a,\,\,u_t\bigl|_{t=0}=b,
    \end{array}
    \right.
\end{equation}
where the functions $a=a(x)$, $b=b(x)$ are the initial conditions
describing spatial distributions of the absorption of laser
radiation energy and of its temporal change.  Also we assume that
there is a set of transducers  running over a surface $S$
measuring the pressure wave and collecting data. A transducer at
at point $x$ at time $t$  measures the value $u(x,t)$. Then we can
suppose that
\begin{equation}
\label{observation} F:=u\bigl|_{S\times[0,T]}.
\end{equation}
is known for some $T>0$. The inverse problem (IP) is to find the
initial data $a$, $b$, using data $F(x,t)$ measured by
transducers. In our setting we assume that the surface $S$ is a
unit sphere and $\supp{a},\supp{b}$ are located inside $S$.

In the present paper we continue the the work started in
\cite{BLMM}. Our goal here is to develop new algorithms that
improve the speed of calculations and the clarity of the resulting
image. In \cite{BLMM} we presented such algorithms for three
($n=3$) and two ($n=2$) dimensional situation. The numerical
realization of these algorithms is quite involved in
multidimensional situation, that is why here we restrict ourselves
for the case $n=1$. On the one hand,  the one-dimensional
situation is rather simple, on the other hand, for numerical
simulation in $n=1$ we preset the simplified version of an
algorithm described in \cite{BLMM}. And we also point out the way
how this simplifications can be made in cases $n=2$ and $n=3$.

In the second section we derive the Fourier representation for the
forward problem in the cases $n=1,2,3.$ In the third section we
describe the algorithm of solving IP in one-dimensional situation.
In the forth section the modifications needed for simplifications
of the algorithm in multidimensional cases are described. In the
last section we demonstrate the results of numerical simulation in
the case $n=1.$

\section{Solving of Forward Problem.}

Since $\supp a,b\subset B_1=\{x\in \mathbb{R}^n\,|\,|x|<1\}$ and
the speed of sound in (\ref{wave_eq}) equal to one, we have that
\begin{equation}\label{cil}
u(x,t)=0,\quad |x|=T+1,\,t\in (0,T).
\end{equation}
We use the separation of variables (Fourier) method to solve the
initial boundary value problem (\ref{wave_eq}),(\ref{cil}), we
look for the solution in a form
\begin{equation}\label{f1}
u(x,t)=\sum_{k=1}^{\infty}X_k(x)T_k(t),
\end{equation}
where $X_k(x)$ are solution to the following spectral problem
\begin{equation}\label{f2}
    \left\{\begin{array}l -\Delta X_k =\lambda_k X_k,\quad |x|<1+T,\\
X_k\bigl|_{|x|=1+T}=0,
\end{array}
\right.
\end{equation}
and $T_k(t)$ are defined by
\begin{equation}\label{f3}
    \left\{\begin{array}l -T''_k =\lambda_k T_k,\quad 0<t<T,\\
\sum_{k=1}^{\infty}X_k(x)T_k(0)=a(x),\quad |x|<1+T, \\
\sum_{k=1}^{\infty}X_k(x)T'_k(0)=b(x),\quad |x|<1+T.
\end{array}
\right.
\end{equation}
From (\ref{f3}) we have that
\begin{equation}\label{ft}
T_k(t)=a_k\cos\sqrt{\lambda_k}t +
\frac{b_k}{\sqrt{\lambda_k}}\sin\sqrt{\lambda_k}t,
\end{equation}
where
\begin{equation}\label{fc}
a_k=\frac{(a,X_k)_{L_2}}{\|X_k\|_{L_2}},\quad
b_k=\frac{(b,X_k)_{L_2}}{\|X_k\|_{L_2}}.
\end{equation}
Solution to (\ref{f2}) has the following form in different
dimensions:
\begin{equation}\label{d1}
X_k(x)=\sin\Bigl(\frac{x+(1+T)}{2(1+T)}k\pi \Bigr),\quad \lambda_k=
\Bigl(\frac{k\pi}{2(1+T)}\Bigr)^2,\quad n=1,
\end{equation}
\begin{equation}\label{d2}
X_{kl}(x)=J_k\Bigl(\frac{\nu_k^lr}{1+T}\Bigr)e^{ik\varphi},\quad \lambda_{kl}=
\Bigl(\frac{\nu_k^l}{1+T}\Bigr)^2,\quad n=2,
\end{equation}
\begin{equation}\label{d3}
X_{klm}(x)=\frac{1}{\sqrt{r}}J_{k+1/2}\Bigl(\frac{\xi_k^lr}{1+T}\Bigr)
e^{im\varphi}P^m_k(\cos\theta),\quad \lambda_{klm}= \Bigl(\frac{\xi_k^l}{1+T}\Bigr)^2,\quad n=3,
\end{equation}
where $J_k$, $J_{k+1/2}$ -- Bessel functions of the first kind,
$\nu_k^l$ and  $\xi_k^l$ -- their roots, $P^m_k$ -- associated
Legendre polynomials:
$$
J_\alpha(x)=\sum_{m=0}^{\infty}\frac{(-1)^m}{m!\Gamma(m+\alpha+1)}\Bigl(\frac{x}{2}\Bigr)^2,
$$
$$
P^m_k(x)=\frac{(-1)^m}{2^kk!} (1-x^2)^{m/2}\frac{d^{l+k}}{dx^{l+k}}(x^2-1)^k.
$$
Thus the representation (\ref{f1}) gives a solution to
(\ref{wave_eq}), (\ref{cil}) in the region $\left\{(x,t)\,|\,
|x|<T+1,\, t\in (0,T) \right\}$, with $X_k$ defined by
(\ref{d1})-(\ref{d3}) and $T_k$ by (\ref{ft}).

\section{Solving of IP for $n=1$.}\label{sip1}

The IP is to find the initial values $a(x)$, $b(x)$ of
(\ref{wave_eq}) using the observation $F(x,t)$ measured on unit
sphere:
$$
F(x,t)=u(x,t), \quad |x|=1,\ 0<t<T.
$$

In 1d case the solution $u$ in $\left\{(x,t)\,|\, |x|<T+1,\, t\in
(0,T) \right\}$ is given by (\ref{f1}), (\ref{ft}), (\ref{fc}),
(\ref{d1}). The latter yields the following equality
\begin{equation}\label{ip}
\sum_{k=1}^{\infty}X_k(x)(a_k\cos\sqrt{\lambda_k}t +
\frac{b_k}{\sqrt{\lambda_k}}\sin\sqrt{\lambda_k}t)=F(x,t),\quad
|x|=1,\ t<T.
\end{equation}
We would like to determine coefficients $a_k$ and $b_k$ from
(\ref{ip}) and restore functions $a$ and $b$ as Fourier series:
$$
a(x)=\sum_{k} a_k X_k(x),\quad b(x)=\sum_{k} b_k X_k(x).
$$
We restrict ourselves to the case $n=1$ and $b=0$, the numerical
testing in the last section was conducted specifically in this
situation. Other cases require some modifications, we postpone
them for further consideration.

In one dimensional case (\ref{ip}) becomes
\begin{equation}
\label{ip1}
\left\{\begin{array}l
\sum_{k=1}^{\infty}\sin\Bigl(\frac{2+T}{2(1+T)}k\pi \Bigr)(a_k\cos\frac{k\pi}{2(1+T)}t + \frac{b_k}{\sqrt{\lambda_k}}\sin\frac{k\pi}{2(1+T)}t)=F(1,t),\\
\sum_{k=1}^{\infty}\sin\Bigl(\frac{T}{2(1+T)}k\pi \Bigr)(a_k\cos\frac{k\pi}{2(1+T)}t + \frac{b_k}{\sqrt{\lambda_k}}\sin\frac{k\pi}{2(1+T)}t)=F(-1,t).
\end{array}
\right.
\end{equation}

Using the equality
$\frac{T}{2(1+T)}k\pi=k\pi-\frac{2+T}{2(1+T)}k\pi$ we get what
$$
F(1,t)=u_1(t)+u_2(t),\quad F(-1,t)=-u_1(2(T+1)-t)+u_2(2(T+1)-t),
$$
where
$$
u_1(t)=\sum_{k=1}^{\infty}a_k\sin\Bigl(\frac{2+T}{2(1+T)}k\pi
\Bigr)\cos\frac{k\pi}{2(1+T)}t,
$$
$$
u_2(t)=\sum_{k=1}^{\infty}b_k\sin\Bigl(\frac{2+T}{2(1+T)}k\pi
\Bigr)\sin\frac{k\pi}{2(1+T)}t.
$$
On assuming that $b=0$ and introducing the notation
$$
\hat F(t)= \left\{\begin{array}l
F(1,t),\quad 0<t<1+T,\\
F(-1,2(1+T)-t),\quad 1+T<t<2(1+T),
\end{array}
\right.
$$
we then arrive at the following relation
\begin{equation}\label{ip11}
\hat F(t)=\sum_{k=1}^{\infty}a_k\sin\Bigl(\frac{2+T}{2(1+T)}k\pi \Bigr)\cos\frac{k\pi}{2(1+T)}t.
\end{equation}
The left hand side of the last equatity is known, therefore  we
can find $a_k$ as a Fourier coefficients of the basis
$\left\{\cos\frac{k\pi}{2(T+1)}t\right\}_{k=0}^\infty$ in
$L_2(0,2(T+1))$: .
\begin{equation}\label{s1}
a_k=\frac{\hat F_k}{\sin\Bigl(\frac{2+T}{2(1+T)}k\pi
\Bigr)},\quad\text{where } \hat F_k=\frac{(\hat
F,\cos\frac{k\pi}{2(1+T)}t )_{L_2(0,2(T+1))}}{T+1}.
\end{equation}
From (\ref{s1}) we can easily see that if, we take $T=2$ then the
denominator in the first formula vanishes for $k=3,6,9,\ldots$.
Therefore we could not use (\ref{s1}) to reconstruct $a$. It is
also the case if $k=(T+1)n$ for some $n\in \mathbb{N}$. To improve
this we need the following lemma:
\begin{lemma}
\label{Lem1} Let $T\in \mathbb{N}$ be fixed, the function $a$ with
$\supp{a}\subset (-1,1)$ admits the expansion
$$
a(x)=\sum_{n=1}^{\infty} a_n\sin\Bigl(\frac{x+1}{2} n\pi\Bigr).
$$
Define $\tilde a(x)$ by the rule:
    $$
    \tilde a(x)= \left\{\begin{array}l
    a(x),\quad -1<x<1,\\
    0,\quad 1<|x|<T+1,
    \end{array}
    \right.
    $$
so $\tilde a(x)$ admits the representation
$$
\tilde a(x)=\sum_{k=1}^{\infty} \tilde
a_k\sin\Bigl(\frac{x+1+T}{2(T+1)} k\pi\Bigr).
$$
Then
\begin{equation}\label{s2}
a(x)=\left\{\begin{array}l
\frac{T+1}{T} A(x),\quad\text{$T$ is even}\\
A(x),\quad\text{$T$ is odd},
\end{array}
\right.
\end{equation}
where
\begin{equation}
\label{A_func}
A(x)= \sum_{k\in \mathbb{N},\,k\not= (T+1)n} \tilde
a_k\sin\Bigl(\frac{x+1+T}{2(T+1)} k\pi\Bigr).
\end{equation}
\end{lemma}
\begin{proof}
We express $\tilde a_k$ via $a_n$:
\begin{eqnarray*}
\tilde a_k=\frac{1}{T+1}\int\limits_{-(T+1)}^{T+1} \tilde a(x)\sin\Bigl(\frac{x+1+T}{2(T+1)} k\pi\Bigr)\,dx= \\
=\frac{1}{T+1}\sum_{n=1}^{\infty}\int\limits_{-1}^{1}  a_n\sin\Bigl(\frac{x+1}{2} n\pi\Bigr)\sin\Bigl(\frac{x+1+T}{2(T+1)} k\pi\Bigr)\,dx=\\
= \left\{\begin{array}l \frac{2}{\pi}\sin\Bigl(\frac{k\pi T}{2(T+1)} \Bigr) \sum\limits_{n=1}^{\infty} \frac{n a_n}{(T+1)^2n^2-k^2},\quad\text{ $n+k$ is even,} \\
0,\quad\text{$n+k$ is odd}.
\end{array}
\right.
\end{eqnarray*}
The last formula works if only $(T+1)n\not=k$ for any $n$, if for
some $n_0$ we have $(T+1)n_0=k$, then
$$
\tilde a_k =\frac{1}{T+1}a_{n_0}\cos\Bigl(\frac{Tn_0}{2}\pi\Bigr).
$$
From the last equality we derive that
\begin{eqnarray*}
\tilde a(x)=A(x)+ \frac{1}{T+1}\sum_{n=1}^{\infty} a_n \cos\Bigl(\frac{Tn}{2}\pi\Bigr)\sin\Bigl(\frac{x+1+T}{2} n\pi\Bigr)=\\
A(x)+\frac{1}{T+1}\sum_{n=1}^{\infty}\frac{a_n}{2} \Bigl(\sin\Bigl(\frac{x+1}{2} n\pi\Bigr)+\sin\Bigl(\frac{x+1}{2} n\pi+Tn\pi\Bigr)\Bigr)=\\
=A(x)+\frac{1}{(T+1)2}\Bigl( a(x)+(-1)^Ta(x)\Bigr).
\end{eqnarray*}
Taking here $|x|<1$ yields (\ref{s2}), which completes the proof.

\end{proof}

Applying this result to our situation, we can make the following
\begin{proposition}
\label{Prop1} Let $T\in \mathbb{N}$ be fixed. The unknown function
$a$ can be recovered by formula (\ref{s2}), where coefficients
$\tilde a_k$ in representation (\ref{A_func}) are given by
(\ref{s1}).

\end{proposition}

\section{Remarks on solving of IP for $n=2,3$.}

We make use of  (\ref{d2}),(\ref{ft}), (\ref{f1}) to write down an
Fourier expansion of observation (\ref{observation}) in 2d case:
$$
\sum_{k,\,l}J_k\Bigl(\frac{\nu_k^l}{1+T}\Bigr)e^{ik\varphi}
\Bigl(a_{kl}\cos\Bigl(\frac{\nu_k^l}{1+T}t\Bigr) +
b_{kl}\frac{1+T}{\nu_k^l}\sin\Bigl(\frac{\nu_k^l}{1+T}t\Bigr)\Bigr)
= F(t,\varphi),
$$
where $F(t,\varphi)$ is a given  boundary measurements and
$a_{kl}$ and $b_{kl}$ are Fourier coefficients subjected to
determination. Multiplying in $L_2(0,2\pi)$ by
$\frac{1}{\sqrt{2\pi}}e^{ik\varphi}$ we obtain:
\begin{equation}\label{ip2}
\sum_{l}J_k\Bigl(\frac{\nu_k^l}{1+T}\Bigr)
\Bigl(a_{kl}\cos\Bigl(\frac{\nu_k^l}{1+T}t\Bigr) +
b_{kl}\frac{1+T}{\nu_k^l}\sin\Bigl(\frac{\nu_k^l}{1+T}t\Bigr)\Bigr)
= F_k(t),
\end{equation}
where $F_k(t)$ is a Fourier coefficient of $F(t,\varphi)$ w.r.t.
family $\left\{\frac{1}{\sqrt{2\pi}}e^{ik\varphi}\right\}_{k\in
\mathbb{Z}}$ in $L_2(0,2\pi)$.

Equation (\ref{ip2}) is a 2d analog of (\ref{ip11}), but the
arguments of cosine and sine functions are different: instead of
$\frac{k\pi}{2(T+1)}$ we have $\frac{\nu_k^l}{T+1}$, where
$\nu_k^l$ is a root of Bessel function $J_k$.

On observing that the system
$\left\{\cos\Bigl(\frac{k\pi}{2(T+1)}t\Bigr)\right\}_{k=0}^\infty$
is orthogonal in $L_2(0,2(T+1))$, we derived formula for $a_k$
(\ref{s1}) from (\ref{ip11}). But now the system
$\left\{\cos\Bigl(\frac{\nu_k^l}{T+1}t\Bigr)\right\}_{l=1}^\infty$
is not a basis. At the same time this family is not ''very bad".


We remind that
$$
J_k(x)=\sqrt{\frac{2}{\pi
x}}\cos\Bigl(x-\frac{k\pi}{2}-\frac{\pi}{4}\Bigr)+\OB\limits_{x\to\infty}\Bigl(\frac{1}{x^{3/2}}\Bigr),
$$
and therefore
$$
\nu_k^l\approx
\frac{k\pi}{2}+\frac{\pi}{4}+\frac{2l-1}{2}\pi,\quad l=1,2,\ldots,
$$
which means that
$$
\nu_k^l- \nu_k^{l+1} \to \pi\quad\text{if $l\to\infty$}.
$$
All aforesaid and \cite{AI}, (see Chapter 4, Theorem II.4.1)
allows us to conclude that the system
$\Bigl\{\cos\Bigl(\frac{\nu_k^l}{T+1}t\Bigr)\Bigr\}_{l=1}^\infty$
is minimal in $L_2(0,T+1)$ and  system
$\Bigl\{\cos\Bigl(\frac{\nu_k^l}{T+1}t\Bigr),\
\sin\Bigl(\frac{\nu_k^l}{T+1}t\Bigr)\Bigr\}_{l=1}^\infty$ is
minimal in $L_2(-T-1,T+1)$. If $b=0$ then
\begin{equation}\label{s3}
a_{kl}=\frac{ F_{kl}}{J_k\Bigl(\frac{\nu_k^l}{1+T} \Bigr)},\quad
F_{kl}=\frac{\Bigl( F_k(t),u_{kl}(t)
\Bigr)_{L_2(0,T+1)}}{\|u_{kl}(t) \|_{L_2(0,T+1)}},
\end{equation}
where the system $\{u_{kl}(t)\}_{l=1}^{\infty}$ is bi-orthogonal
to
$\Bigl\{\cos\Bigl(\frac{\nu_k^l}{T+1}t\Bigr)\Bigr\}_{l=1}^\infty$
in $L_2(0,T+1)$. Formula (\ref{s3}) is 2d analog of (\ref{s1}),
but instead of $\sin\Bigl(\frac{2+T}{2(1+T)}k\pi \Bigr)$ in
denominator we have $J_k\Bigl(\frac{\nu_k^l}{1+T} \Bigr)$. It has
an advantage over one-dimensional case because the denominator of
(\ref{s3}) does not vanish for integer values of $T$. Disadvantage
is that denominator tend to zero while $l\to\infty$ which makes
calculation  sensitive to accuracy.

Similar reasoning can be used in a three-dimensional situation,
but instead of (\ref{d2}) we use (\ref{d3}).


\section{Numerical experiment}

In this section we provide the results on numerical testing for
the one-dimensional problem. Namely, we consider (\ref{wave_eq})
with $n=1$, $b=0$, take $T=2$ in (\ref{cil}). First we generate
data for IP i.e. observation (\ref{observation}), by solving
forward problem with the use of (\ref{f1}), (\ref{ft}),
(\ref{fc}), (\ref{d1}).

After solving forward problem we discretized obtained data
(function $F(\pm1,t)$) to get closer to the case where the data is
the result of a real experiment, i.e. we replace $F(\pm1,t)$ by
$F(\pm1,\frac{k}{N})$, $k=1,\ldots,2N$ for some large $N$, we also
add some noise to the data to be more realistic.
\begin{figure}[h]
    \begin{center}
        \begin{minipage}[h]{0.49\linewidth}
            \includegraphics[width=1\linewidth]{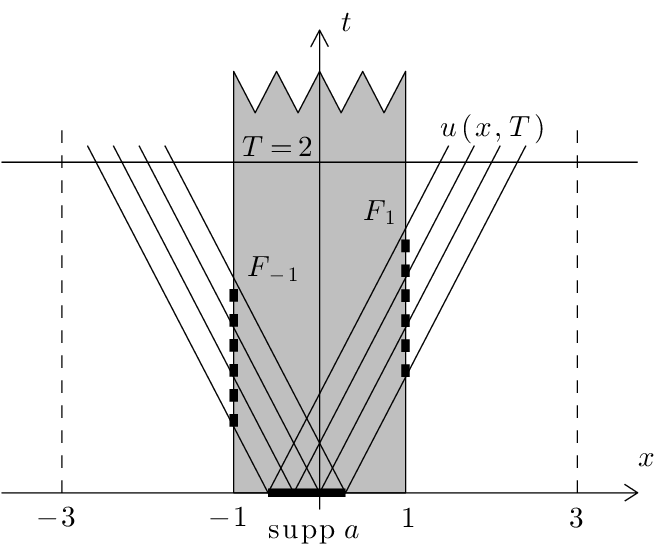}
            \caption{Space-time cylinder.} 
            \label{ris-st} 
        \end{minipage}
        \hfill
        \begin{minipage}[h]{0.5\linewidth}
            \includegraphics[width=1\linewidth]{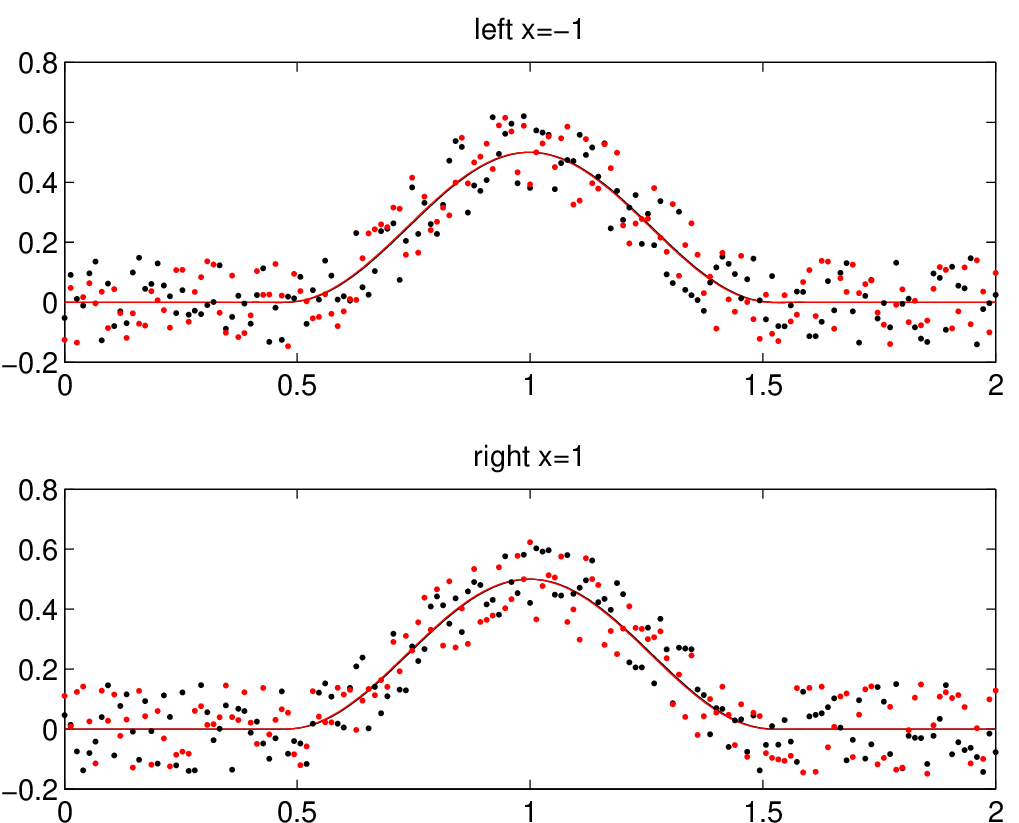}
            \caption{Observation.}
            \label{ris-ob}
        \end{minipage}
    \end{center}
\end{figure}

For $a(x)=1/2+1/2\cos(2\pi x)$ observation $F_{\pm1}$ is
represented on a picture. Before applying our algorithm, we
checked the possibility of solving this problem by brute force,
i.e. minimizing the least-square-error for a Fourier polynomial.
As it is shown on the picture, the result was negative:
\begin{figure}[h!]
    \begin{center}
        \begin{minipage}[h]{0.49\linewidth}
            \includegraphics[width=1\linewidth]{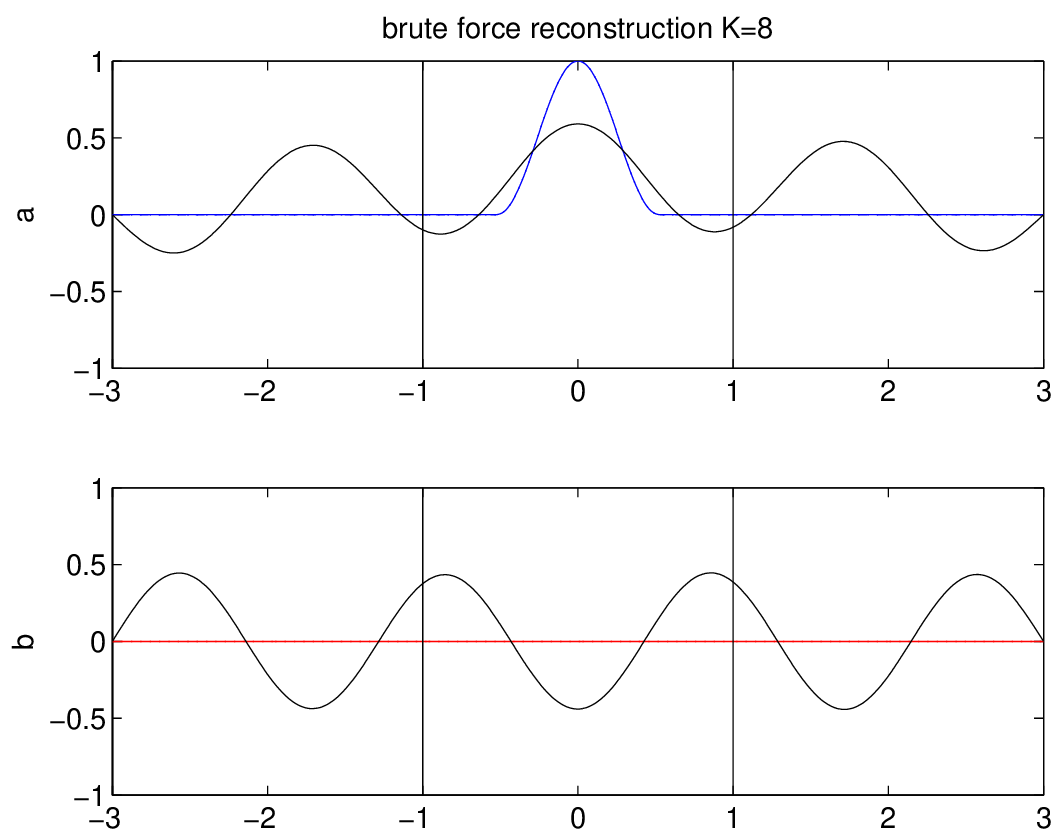}
        \end{minipage}
        \hfill
        \begin{minipage}[h]{0.5\linewidth}
            \includegraphics[width=1\linewidth]{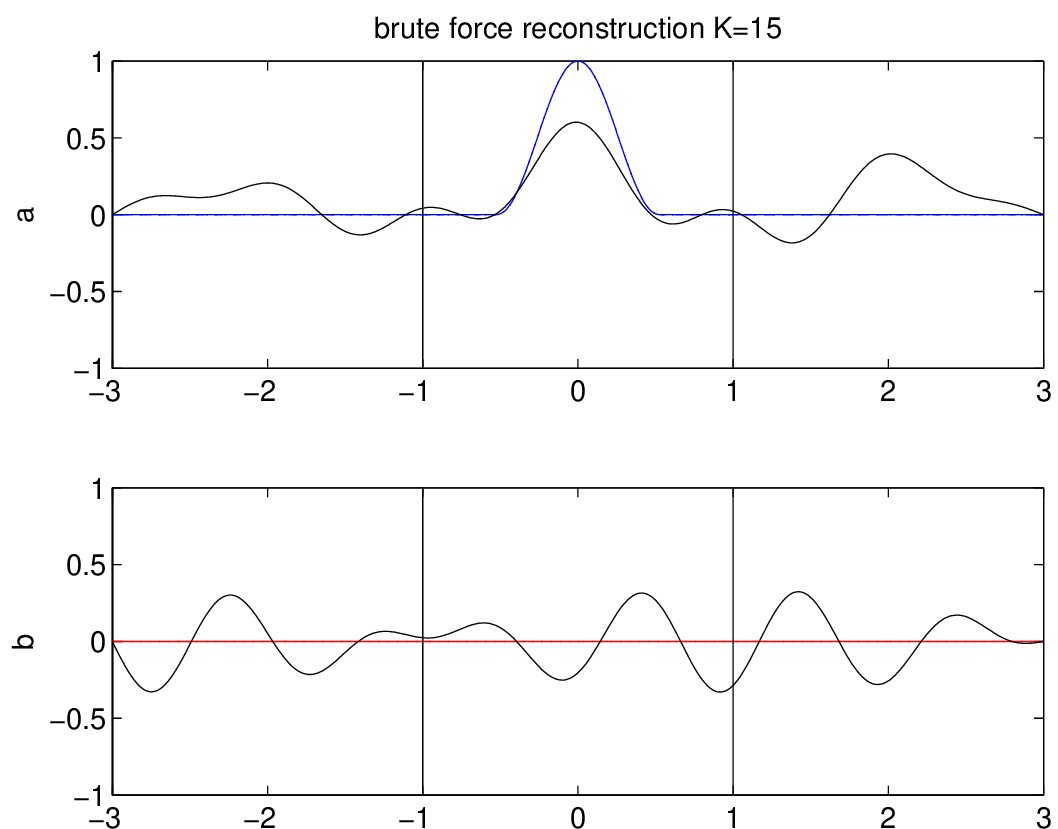}
        \end{minipage}
    \end{center}\caption{Minimizing least-square-error.} 
\label{ris-bi}
\end{figure}

We applied our algorithm, described in section \ref{sip1} with
$T=2$. As was mentioned, the reconstruction of Fourier
coefficients $a_k$ using (\ref{s1}) is possible only for
$k\not=3n$ where $n=1,2,\ldots$ because for $k=3n$ the denominator
in (\ref{s1}) vanishes. Making use of Lemma \ref{Lem1} and
Proposition \ref{Prop1} we get the following answer:
$$
a(x)=\frac{3}{2} A(x),\quad\text{where}
\quad A(x)= \sum_{k\not= 3n}  a_k\sin\Bigl(\frac{x+3}{6} k\pi\Bigr).
$$

Reconstructed data   in the 'error' case  with a 50 Fourier
coefficients are represented on figures
\begin{figure}[h!]
    \begin{center}
        \begin{minipage}[h]{0.49\linewidth}
            \includegraphics[width=1\linewidth]{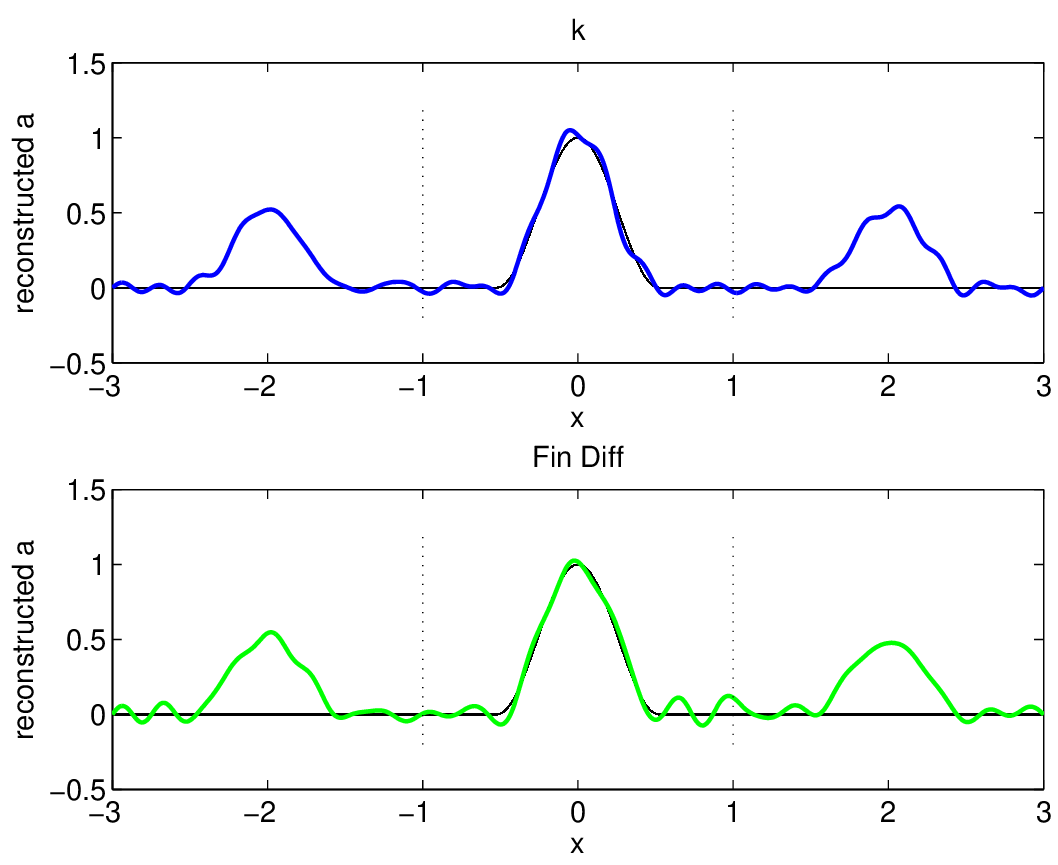}
                        \caption{Smooth wave.} 
                        \label{ris-sin} 
        \end{minipage}
        \hfill
        \begin{minipage}[h]{0.5\linewidth}
            \includegraphics[width=1\linewidth]{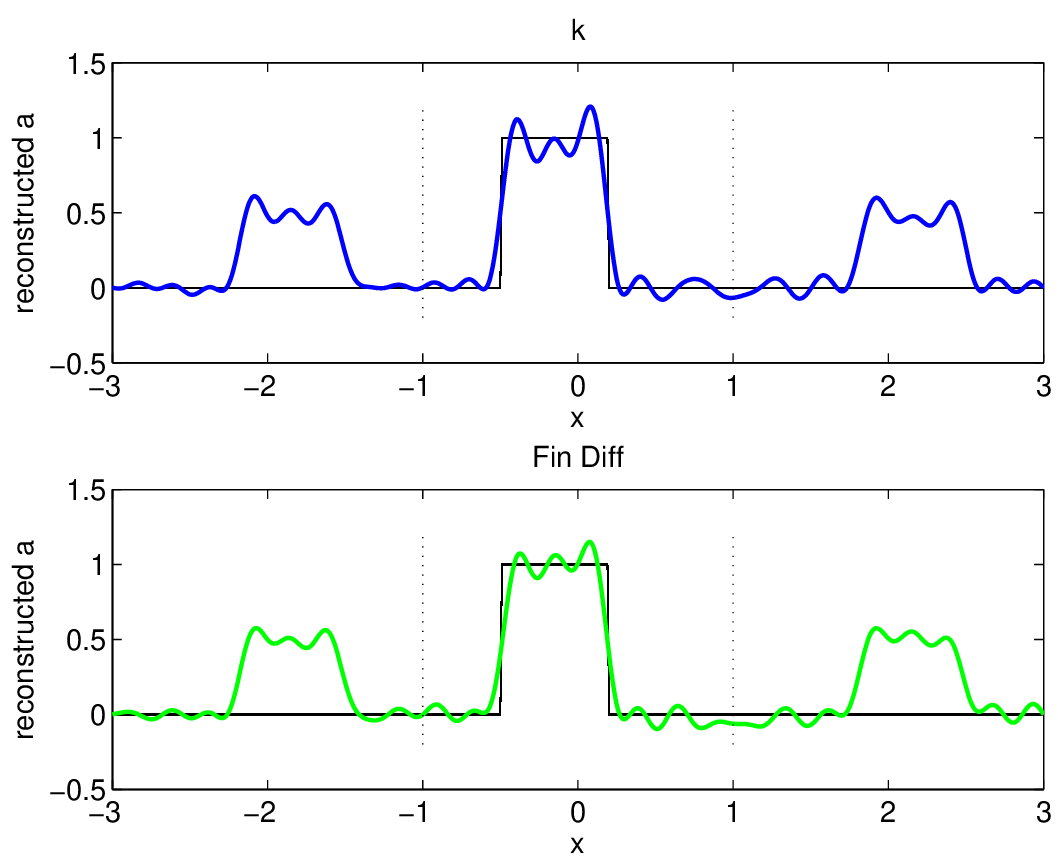}
                        \caption{Step function.}
                        \label{ris-step}
        \end{minipage}
    \end{center}
\end{figure}

The last two figures demonstrate comparison of our algorithm
(first row) with the reconstruction by solving wave equation
backward in time by finite-difference scheme (second row). On the
right figure we reconstruct initial condition $a(x)$ having
nonsmooth profile.

\noindent{\bf Acknowledgments}

Alexandr Mikhaylov and Victor Mikhaylov were supported by RFBR
18-01-00269 and by the Ministry of Education and Science of
Republic of Kazakhstan under grant AP05136197. We thank the
Volkswagen Foundation (VolkswagenStiftung)  program ``Modeling,
Analysis, and Approximation Theory toward application in
tomography and inverse problems'' for kind support and stimulating
our collaboration.

\end{document}